\documentclass[12pt]{article}
\usepackage[utf8]{inputenc}
\usepackage{color}
\usepackage{amsmath,amssymb,amsthm,amsfonts}
\usepackage{float}
\usepackage{graphicx}
\RequirePackage{hyperref}

\newtheorem{remark}{Remark}

\newtheorem{definition}{Definition}
\newtheorem{theorem}{Theorem}

\usepackage{geometry}
\begin{document}
\begin{center}
\textsc{{ \Large Harvest management problem in a fractional logistic equation }\\[0pt]}
\vspace{0.4cm} \hspace{0.2cm} Melani Barrios$^{1,2}$ \hspace{0.5cm} Gabriela Reyero$^{1}$ \hspace{0.5cm} Mabel Tidball$^{3}$ \\[0pt]
\end{center}

\scriptsize                                                                               
 $\,^{1}$  Departamento de Matem\'atica, Facultad de Ciencias Exactas, Ingeniería y Agrimensura, Universidad Nacional de Rosario, Avda. Pellegrini $250$, S$2000$BTP Rosario, Argentina.
 
$\,^{2}$ CONICET, Departamento de Matem\'atica, Facultad de Ciencias Exactas, Ingeniería y Agrimensura, Universidad Nacional de Rosario, Avda. Pellegrini $250$, S$2000$BTP Rosario, Argentina.

$\,^{3}$ CEE-M, Universidad de Montpellier, CNRS, INRA, SupAgro, Montpellier, France.\\

\normalsize

Correspondence should be addressed to melani@fceia.unr.edu.ar\\

\begin{abstract} 
In this article, we study a fractional control problem that models the maximization of the profits obtained by exploiting a certain resource whose dynamics are governed by the fractional logistic equation. Due to the singularity of this problem, we develop different resolution techniques, both for the classical case and for the fractional case. In the last section we perform several numerical simulations to make a comparison between both cases.
\end{abstract}




\section{Introduction}
The logistic equation describes the population growth. The model is initially published by Pierre Verhulst in 1838 \cite{Cu}. The continuous Logistic model is described by a first order ordinary differential equation. The model describes the population growth that may be limited by certain factors like population density \cite{AlSaYo, Au, Pa}. 
The continuous form of the logistic equation is expressed in the form of a nonlinear ordinary differential equation,
\[\dot{x}(t)=rx(t)\left(1-\frac{x(t)}{K}\right).\]
In the above equation, $x(t)$ indicates population at time $t$, $r>0$ represents the Malthusian parameter expressing growth rate of species and $K$ denotes carrying capacity.

Motivated by its applications in different scientific areas (electricity, magnetism, mechanics, fluid dynamics, medicine, etc. \cite{Alm, BaRe, Hil, Kil}), fractional calculus is in development, which has led to great growth in its study in recent decades. The fractional derivative is a nonlocal operator \cite{Die, Pod}, making fractional differential equations good candidates for modeling situations in which is important to consider the history of the phenomenon studied \cite{FeSa}, unlike the models with classical derivative where this is not taken into account. There are several definitions of fractional derivatives. The most commonly used are the Riemann-Liouville fractional derivative and the Caputo fractional derivative. It is important to note that while the Riemann-Liouville fractional derivative \cite{Old}, is historically the most studied approach to fractional calculus, the Caputo fractional derivative is more popular among physicists and scientists due to the fact that the formulation of initial value problems with this type of derivative is more similar to the formulation with classical derivative.

The fractional order logistic equation has been discussed in the literature \cite{BhDaGe, MoQa}. A detailed study of existence, uniqueness, stability and approximate solutions of this equation can be found in \cite{AmNuAnSu, BaReTid, EsEmEs, KSQB, SwKhMa, SyMaSh}.

When the logistic equation is used to describe the natural evolution of a species, it is logical to think about the exploitation of this resource. For this reason we have decided to make a study on the maximization of the exploitation of a certain resource by using fractional derivatives.

The principal objective of this work is to study a fractional control problem that models the maximization of the profits obtained by exploiting a certain resource. The structure of this article is as follows: in section 2, the classical control problem and its solution are presented. In section 3, the fundamental concept of  fractional derivatives, fractional control and variational problems, and the solution of a fractional control problem. The comparison between both problems and numerical approximations of the solutions are discussed in section 4. Finally, the section 5, is dedicated to the conclusions.

\section{A simple optimal control problem}\label{sec2}
A simple classical optimal control model will be presented, which consists of maximizing the extraction or harvest of a certain renewable resource, where the first-order integer derivative appears in the dynamic equation,
\begin{equation}
\left\{
\begin{array}{l}
\max \, \int \limits_{0}^{T} e^{-\delta t} h(x(t)) \, dt\\
\dot{x}(t)=rx(t)\left(1-\frac{x(t)}{K}\right)-h(x(t))\\
x(0)=x_0\\
x(T)=x_{T}\\
h_{min}\leq h(x(t)) \leq h_{max},\\
\end{array}\right.
\label{pclasico}
\end{equation}
where  $T$ represents the final time, $x_0$ the initial condition, $x_T$ the final condition, $h_{min}$ and $h_{max}$ preset minimum and maximum harvest and $e^{-\delta t}$ represents a discount factor with $\delta\geq 0$ the instantaneous annual rate of discount which can be zero.

In this problem the following dynamic equation appears
\begin{equation}
\dot{x}(t)=rx(t)\left(1-\frac{x(t)}{K}\right)-h(t),
\label{eclogisticacosecha}
\end{equation}
which is the logistical growth of a certain resource, where $r>0$ is the intrinsic growth rate, $K>0$ the carrying capacity of the resource and $h(t)$ is the harvest.  
\begin{remark}
In general, the function to be optimized is the following
$$  \left(p -  \frac{c}{x} \right) h(x),$$
where $p$ is the price and $\frac{c}{x}$ the extraction cost, \cite{Cl}. In this work a simplification of this model will be done, since the resolution of this problem has a high level of complexity in its fractional version. 
\end{remark}
Now we will proceed to solve the problem (\ref{pclasico}).

Clearing $h(x(t))$ from the dynamic equation (\ref{eclogisticacosecha}) and replacing it in the functional to maximize results,  
$$\left\{
\begin{array}{l}
\max \, \int \limits_{0}^{T} e^{-\delta t} \left[rx(t)\left(1-\frac{x(t)}{K}\right)-\dot{x}(t)\right] \, dt\\
x(0)=x_0\\
x(T)=x_{T}\\
h_{min}\leq h(x(t)) \leq h_{max},\\
\end{array}\right.
$$
which means, it becomes a variational problem with Lagrangian  
$$L(t,x,\dot{x})=e^{-\delta t} \left[rx(t)\left(1-\frac{x(t)}{K}\right)-\dot{x}(t)\right].$$
Then, its Euler-Lagrange equation results 
$$
\begin{array}{l}
\dfrac{\partial L}{\partial x}-\dfrac{d}{dt} \dfrac{\partial L}{ \partial \dot{x}}=0,
\end{array}
$$
from which is obtained
$$
\begin{array}{l}
r\left(1-\frac{2}{K}x(t)\right)-\delta=0.
\end{array}
$$
That is, the optimal population is the constant 
\begin{equation}
x^{*}(t)=\frac{K}{2}\left(1-\frac{\delta}{r}\right),
\label{condesc1}
\end{equation}
while the optimal harvest is
\begin{equation}
h^{*}(t)=rx^{*}(t)\left(1-\frac{x^{*}(t)}{K}\right),
\end{equation}
then, 
\begin{equation}
h^{*}(t)=\frac{K}{4} \frac{r^2-\delta^2}{r}.
\label{condesch1}
\end{equation}

\begin{remark}
It can be seen that being the solution $x^{*}(t)$ a constant function, it will not be able to verify the imposed boundary conditions. For this reason, it can be said that the problem is singular.

For its resolution, the Nearest Feasible Paths theorem will be used, and its proof can be seen in \cite{Set}.    
\end{remark}

\begin{theorem}\textbf{Nearest Feasible Paths}

Consider the following optimal control problem with one state and one control variable
$$\left\{
\begin{array}{l}
\max \limits_{u} \, \int \limits_{0}^{T} e^{-\delta t} F(x,u) \, dt\\
\dot{x}(t)=f(x,u)\\
x(0)=x_0\\
x(T)=x_T\\
u \in [0,Q]\subset \mathbb{R},
\end{array}\right.
$$
we define the Nearest Feasible Paths $x^{*}(t)$ to a given feasible steady state $(\overline{x},\overline{u})$ as the feasible path starting from $x_0$ in such a way that
\[\left|x^{*}-\overline{x}\right|\leq \left|x(t)-\overline{x}\right| \, \, \forall t>0,\]
for all feasible paths $x(t)$ starting from $x_0$. Further, we define the nearest approach segment $\hat{x}^{*}(t), \, t \in [0, t_{min}]$, as the feasible segment starting from $x_0$ in such a way that
\[\left|\hat{x}^{*}-\overline{x}\right|\leq \left|x(t)-\overline{x}\right| \, \, \forall t\in [0, t_{min}],\]
for all feasible segments $x(t)$ starting from $x_0$ and where $t_{min}$ is defined as the first time $\hat{x}^{*}(t)$ reaches $\overline{x}$. If $\overline{x}$ were the optimal stationary equilibrium and if the optimal path for the problem is the Nearest Feasible Paths, then the nearest approach segment from $x_0$ to $\overline{x}$, followed by a stay at $\overline{x}$, followed by an exit from $\overline{x}$ as late as possible ($t_{max}$) to attain $x_T$ at time $T$ is optimal.
\end{theorem}

 With this result, the method consists in finding $x_{a}(t)$ the solution to the problem that begins in $x(0)=x_0$ and uses the minimum harvest $h_{min}$ if $x_0<x^{*}(t)$ or uses the maximum harvest $h_{max}$ if $x_0>x^{*}(t)$, until it reaches the value $x^{*}(t)$ in the time $t_{min}$ to determinate. We call its solution $x_a(t)$ and $t_{min}$ will be the first moment when $x_{a}(t)$ gets the value $x^{*}(t)$.
 
Then it must be found $x_{b}(t)$ solution to the problem with final condition $x (T)=x_T$, that uses the minimum harvest $h_{min}$ if $x_T>x^{*}(t)$ or uses the maximum harvest $h_{max}$ if $x_T<x^{*}(t)$.  The time $t_{max}$ is determined so that it is the first time when $x_b(t)$ gets the value $x^{*}(t)$.

In summary, the following optimal solution is obtained
\begin{equation}
x^{*}(t)= \left\{ \begin{array}{lcc}
             x_{a}(t) &   \mathrm{if}  & 0\leq t \leq  t_{min}   \\
             \\ \frac{K}{2}\left(1-\frac{\delta}{r}\right) &  \mathrm{if} & t_{min} \leq t\leq t_{max} \\
             \\ x_{b}(t) &  \mathrm{if}  & t_{max}  \leq t\leq T
             \end{array}
   \right.
	\label{xclasic}
\end{equation}
and the optimal harvest obtained with this procedure is
\begin{equation}
h^{*}(t)= \left\{ \begin{array}{lcc}
             h_{min} &   \mathrm{if}  & x(t)> x^{*}(t) \\
             \\ \frac{K}{4} \frac{r^2-\delta^2}{r} &  \mathrm{if} & x(t) = x^*(t)\\
             \\ h_{max} &  \mathrm{if}  & x(t)< x^{*}(t),
             \end{array}
   \right.
\label{hclasic}
\end{equation}
taking into account that for it to be a feasible solution it must be satisfied 
\[h_{min}\leq \frac{K}{4} \frac{r^2-\delta^2}{r} \leq h_{max}.\] 
\begin{remark}
The problem without discount is obtained as a particular case of the given, writing $\delta=0$ in the equations (\ref{xclasic}) and (\ref{hclasic}). 
\end{remark}

Following, this problem will be solved but in its fractional version, for this we will need to recognize tools of the fractional calculus, which allow us to change from fractional control problems to fractional variational problems and their corresponding fractional Euler-Lagrange equation.

\section{A fractional optimal control problem} \label{sec:frac}
\subsection{Introduction to fractional calculus}
 In this section certain definitions and properties of the fractional calculus will be presented. For more details refer to \cite{Die, Old, Pod}.
\begin{definition}
The Mittag Leffler function with parameters $\alpha , \, \beta$, is defined by 
\begin{equation}
 E_{\alpha,\beta} (z) = \displaystyle \sum_{k=0}^{\infty} \frac{z^k}{\Gamma(\alpha k + \beta)},
\label{mitag}
\end{equation}
for all $z\in \mathbb{C}$.
\end{definition}

\begin{definition}
The Gamma function, $\Gamma: (0, \infty)\rightarrow \mathbb{R}$, is defined by
\begin{equation}
 \Gamma(t) = \int_{0}^{\infty} s^{t-1} e^{-s} \, ds.
\label{gamma}
\end{equation}
\end{definition}
\begin{definition}
The Riemann-Liouville fractional integral operator of order $\alpha \in \mathbb{R}^{+}_{0}$ is defined in $L^1[a,b]$ by
\begin{equation}
 \,_{a}I_{t}^{\alpha} [f] (t) = \dfrac{1}{\Gamma(\alpha)} \int_{a}^{t} (t-s)^{\alpha -1} f(s) \, ds.
\label{frac1}
\end{equation}
\end{definition}
\begin{definition}(Left and Right Riemann-Liouville Fractional Derivatives) \label{defRL} 

The left and right Riemann-Liouville fractional derivatives of order $\alpha \in \mathbb{R}^{+}_{0}$ are defined, respectively, by
\[ \,_{a}D_{t}^{\alpha}[f](t)= \dfrac{1}{\Gamma(n-\alpha)}\dfrac{d^{n}}{dt^{n}}\int_a^{t}(t-s)^{n-1-\alpha}f(s)ds\]
and
\[ \,_{t}D_{b}^{\alpha}[f](t)= \dfrac{(-1)^{n}}{\Gamma(n-\alpha)}\dfrac{d^{n}}{dt^{n}}\int_t^{b}(s-t)^{n-1-\alpha}f(s)ds,\]
with $n=\left\lceil \alpha \right\rceil$, that is if $f \in L^1[a,b]$.
\end{definition}
\begin{definition}(Left and Right Caputo Fractional Derivatives)

The left and right Caputo fractional derivatives of order $\alpha \in \mathbb{R}^{+}_{0}$ are defined, respectively, by
\[ \,_{a}^{C}D_{t}^{\alpha}[f](t)= \dfrac{1}{\Gamma(n-\alpha)}\int_a^{t}(t-s)^{n-1-\alpha}\dfrac{d^{n}}{ds^{n}}f(s)ds\]
and
\[ \,_{t}^{C}D_{b}^{\alpha}[f](t)= \dfrac{(-1)^{n}}{\Gamma(n-\alpha)}\int_t^{b}(s-t)^{n-1-\alpha}\dfrac{d^{n}}{ds^{n}}f(s)ds,\]
with $n=\left\lceil \alpha \right\rceil$, that is if  $\dfrac{d^{n}f}{dt^{n}} \in L^1[a,b]$.
\end{definition}
\begin{remark}An important difference between Riemann-Liouville derivatives and Caputo derivatives is that, being K an arbitrary constant, is
\[ \,_{a}^{C}D_{t}^{\alpha} K= 0,\hspace{1cm} \,_{t}^{C}D_{b}^{\alpha}K=0, \]
however
\[ \,_{a}D_{t}^{\alpha}K=\dfrac{K}{\Gamma (1- \alpha)}(t-a)^{-\alpha}, \hspace{1cm} \,_{t}D_{b}^{\alpha}K=\dfrac{K}{\Gamma (1- \alpha)}(b-t)^{-\alpha}, \]
\[ \,_{a}D_{t}^{\alpha} (t-a)^{\alpha-1}=0, \hspace{1cm} \,_{t}D_{b}^{\alpha}(b-t)^{\alpha-1}=0.\]
In this sense, the Caputo fractional derivatives are similar to the classical derivatives.  
\end{remark}

\subsection{Fractional control and variational problems}
To solve a fractional control problem, tools of fractional variational problems will be used. For this, a brief introduction to them is presented.

 Consider the following problem of the fractional calculus of variations:
find a function $x \in \,_{a}^{\alpha}E$ that optimizes (minimizes or maximizes) the functional
\[J(x)= \int^{b}_{a} L(t,x,\,_{a}^{C}D_{t}^{\alpha}x) \, dt,\]
with a Lagrangian $L \in C^1([a,b]\times \mathbb{R}^2)$ and 
\[\,_{a}^{\alpha}E=\{ x: [a,b] \rightarrow \mathbb{R}: x \in C^1([a,b]), \, \,_{a}^{C}D_{t}^{\alpha}x \in C([a,b]) \},\]
subject to the boundary conditions $x(a)=x_a \, ,\, \, x(b)=x_{b}$.
 
Now the Euler-Lagrange equation for this problem will be stated, its proof is in \cite{LaTo}.

\begin{theorem}
Let $x$ be an optimizer of $J$ in $\,_{a}^{\alpha}E$ with $L\in C^{2}\left([a,b]\times\mathbb{R}^{2}\right)$ subject to boundary conditions $x(a)=x_a \, ,\, \, x(b)=x_{b}$, then $x$ satisfies the fractional Euler-Lagrange differential equation
\begin{equation}
\dfrac{\partial{L}}{\partial x}+\,_{t}^{C}D_{b}^{\alpha}\dfrac{\partial{L}}{\partial \,_{a}^{C}D_{t}^{\alpha}x}=0.
\label{EulerLagrange}
\end{equation}
\label{teoEulerLagrange}
\end{theorem}

\subsection{Fractional model}
Below, the fractional version of the problem that has been proposed in the first section will be shown. This means, it is the same problem but where derivatives of a fractional order appear,
\begin{equation}
\left\{
\begin{array}{l}
\max \, \int \limits_{0}^{T} e^{-\delta t} h(x(t)) \, dt\\
\,^{C}_{0} D^{\alpha}_{t}\left[x\right](t)=rx(t)\left(1-\frac{x(t)}{K}\right)-h(x(t))\\
x(0)=x_0\\
x(T)=x_{T}\\
h_{min}\leq h(x(t)) \leq h_{max}.\\
\end{array}\right.
\label{pfraccionario}
\end{equation}
The only difference with the problem exposed above, is that in the dynamic equation
\begin{equation}
\,^{C}_{0} D^{\alpha}_{t}\left[x\right](t)=rx(t)\left(1-\frac{x(t)}{K}\right)-h(x(t)),\\
\label{dinamicafrac}
\end{equation}
the first-order derivative no longer appears, but now intervenes $\,^{C}_{0} D^{\alpha}_{t}\left[x\right](t)$ the left Caputo fractional derivative of order $0<\alpha\leq 1$.

On one hand, the derivatives with fractional order contains partially or totally the history, temporary future or the spatial behavior of the function, averaged in some way. This transforms the fractional differential equations on suitable candidates for the modeling of memory phenomenon or subsequent effects, those in which what happens at a point on the space or at an instant of time depends on an interval (spatial or temporal) that has the point or the instant.

The Riemann-Liouville fractional derivative had a determining role in the developing of the fractional calculus theory, and was used successfully in strictly mathematical applications. But when it was trying to carry out mathematical modeling of real physical phenomena using fractional differential equations, the problem of the initial conditions also of fractional order emerged. These types of conditions are not physically interpretable and presents a considerable obstacle when making practical use of fractional calculus. The Caputo differential operator, in contrast to the Riemann-Liouville operator, uses derivatives of integer order as initial conditions, that is, initial values that are physically interpretable as in the models with integer derivatives. The definition that follows represented a notable practical advance in the study of physical phenomena such as those of the viscoelastic type and others.

Finally, the fractional derivative at $t$ of a function $x$ is a non-local operator, depending on past values of $x$ (left derivatives) or future values of $x$ (right derivatives). In physics, the right fractional derivative of $x(t)$ is interpreted as a future state of the process $x(t)$. For this reason, the right derivative is usually neglected in applications, when the present state of the process does not depend on the results of the future development. However, the left fractional derivative of $x(t)$ is interpreted as a past state of the process $x(t)$, in which memory effects intervene.

Since the evolution of a certain resource depends on its past, we have decided to choose the left Caputo fractional derivative for modeling its evolution, see \cite{CaTo}.

Now the problem (\ref{pfraccionario}) will be solved.

Clearing $h(x(t))$ from the dynamic equation (\ref{dinamicafrac}) and replacing it in the functional to maximize, it results, 
$$\left\{
\begin{array}{l}
\max \, \int \limits_{0}^{T} e^{-\delta t} \left[rx(t)\left(1-\frac{x(t)}{K}\right)-\,^{C}_{0} D^{\alpha}_{t}\left[x\right](t)\right] \, dt\\
x(0)=x_0\\
x(T)=x_{T}\\
h_{min}\leq h(x(t)) \leq h_{max}.\\
\end{array}\right.
$$
Again it results a variational problem, now fractional, with fractional Lagrangian
$$L(t,x,\,^{C}_{0} D^{\alpha}_{t}\left[x\right])=e^{-\delta t} \left[rx(t)\left(1-\frac{x(t)}{K}\right)-\,^{C}_{0} D^{\alpha}_{t}\left[x\right](t)\right],$$
belonging to $C^{2}\left([0,T]\times\mathbb{R}^{2}\right)$.

Using theorem \ref{teoEulerLagrange}, its fractional Euler-Lagrange equation (\ref{EulerLagrange}) results,
$$
\begin{array}{l}
\dfrac{\partial L}{\partial x}+ \,^{C}_{t} D^{\alpha}_{T}\left[ \dfrac{\partial L}{ \partial \,^{C}_{0} D^{\alpha}_{t}\left[x\right] } \right] =0, 
\end{array}
$$
from which is obtained
$$
\begin{array}{l}
e^{-\delta t} r \left(1-\frac{2}{K}x(t)\right)-\delta (T-t)^{1-\alpha}e^{-\delta T} E_{1,2-\alpha}(\delta (T-t))=0,
\end{array}
$$
where $E_{1,2-\alpha}(\delta (T-t))$ is the Mittag-Leffler function of two parameters.

We can conclude that the optimal population is
\begin{equation}
\begin{array}{l}
x^{*}_{\alpha}(t)=\frac{K}{2}\left(1-\frac{\delta}{r} (T-t)^{1-\alpha} e^{-\delta (T-t)} E_{1,2-\alpha}(\delta (T-t))\right).

\end{array}
\label{condescalfa}
\end{equation} 
And the optimal harvest is obtained from this expression
\begin{equation}
h^{*}_{\alpha}(t)=rx^{*}_{\alpha}(t)\left(1-\frac{x^{*}_{\alpha}(t)}{K}\right)-\,^{C}_{0} D^{\alpha}_{t}\left[x^{*}_{\alpha}\right](t).
\label{condeschalfa}
\end{equation}
Note that in the classical case the optimal population is constant, which is not happening in this case, therefore in the optimal harvest a term appears with the Caputo derivative of $x^{*}_{\alpha}(t)$, that due to its difficulty we must calculate it numerically.
\begin{remark}
We can see that, as in the classical case, the optimal solution $x^{*}_{\alpha}(t)$, although it is not a constant, it will also not verify the established boundary conditions. This means that we are once again faced with a singular problem. 

For its resolution, it would be possible to resort again to a Nearest Feasible Paths theorem, but in a fractional version. However, it has not been carried out yet because the proof of the theorem requires the use of a fractional Green theorem, which at the moment is only available for rectangular regions\cite{OdMaTo}, and its version is necessary for all types of regions. 
\end{remark}
 
\begin{remark}
Here is also the solution of the problem without discount ($\delta=0$) as a particular case of the solution of the problem (\ref{pfraccionario}). 

The optimal population is
\begin{equation}
\begin{array}{l}
x^{*}_{\alpha}(t)=\frac{K}{2}\\
\end{array}
\label{sindescalfa}
\end{equation}
and the optimal harvest is obtained from this expression
\begin{equation}
h^{*}_{\alpha}(t)= \frac{rK}{4}.
\label{sindeschalfa}
\end{equation}
It can be observed that the solution is the same for the classical case, we have that for the case $\delta=0$ the solution and the optimal harvest is the same for all $0<\alpha\leq 1$. The reason  of this is that there is no time-dependent factor that multiplies $\,^{C}_{0} D^{\alpha}_{t}\left[x^{*}_{\alpha}\right](t)$ in the fractional variational problem resulting from taking $\delta=0$, 
$$\left\{
\begin{array}{l}
\max \, \int \limits_{0}^{T} \left[rx(t)\left(1-\frac{x(t)}{K}\right)-\,^{C}_{0} D^{\alpha}_{t}\left[x\right](t)\right] \, dt\\
x(0)=x_0\\
x(T)=x_{T}\\
h_{min}\leq h(x(t)) \leq h_{max},\\
\end{array}\right.
$$
with Lagrangian $L(t,x,\,^{C}_{0} D^{\alpha}_{t}\left[x\right])= rx(t)\left(1-\frac{x(t)}{K}\right)-\,^{C}_{0} D^{\alpha}_{t}\left[x\right](t)$.

Its fractional Euler-Lagrange equation is
$$
\begin{array}{l}
\dfrac{\partial L}{\partial x}+ \,^{C}_{t} D^{\alpha}_{T}\left[ \dfrac{\partial L}{ \partial \,^{C}_{0} D^{\alpha}_{t}\left[x\right] } \right] =0,\\
\end{array}
$$
which means that
$$
\begin{array}{l}
r \left(1-\frac{2}{K}x(t)\right)+\,^{C}_{t} D^{\alpha}_{T}\left[(-1)\right] =0,
\end{array}
$$
then, 
$$
\begin{array}{l}
 r \left(1-\frac{2}{K}x(t)\right)=0,
\end{array}
$$
since the Caputo derivative of a constant is also zero, as in the case of the classical derivative, then the same Euler-Lagrange equation is obtained for all $0<\alpha\leq 1$ and its solution is independent of that value.
\end{remark}

To perform a graphic analysis of the solutions to the discounted problem in both the fractional and classical cases, we must see a specific example. 

\section*{Example and comparison}
Consider $r,\,K,\,\delta,\,x_0,\, x_T$, $T$, $h_{min}$ and $h_{max}$ in given values, as in the example \textit{The Pacific Halibut Fishery}, in \cite{Cl}.

As well, it is considered a fixed value of $\alpha$ to evade the problem of needing a fractional Nearest Feasible Paths theorem, being able to obtain the boundary conditions $x_0$ and $x_{10}$ of this fractional solution, in order to make a comparison.

The following problem will be considered,  
$$\left\{
\begin{array}{l}
\max \, \int \limits_{0}^{10} e^{-0.01 t} h(x(t)) \, dt\\
\,^{C}_{0} D^{\alpha}_{t}\left[x\right](t)=0.71\,x(t)\left(1-\frac{x(t)}{80.5}\right)-h(x(t))\\
x(0)=38.6896\\
x(10)=40.25\\
10 \leq h(x(t)) \leq 15,\\
\end{array}\right.
$$
A comparison of the results will be made between $\alpha=1$ (classical version) and $\alpha=0.6$ (fractional version).

Solving the problem in its classical version as in the previous section, using (\ref{condesc1}) and (\ref{condesch1}) we have
\begin{equation}\label{xestrella}
x^{*}(t)=39.6831,
\end{equation}
\begin{equation}
h^{*}(t)=14.2859.
\end{equation}
It can be observed that $x^{*}(t)$ does not verify established boundary conditions, therefore, the Nearest Feasible Paths theorem will be used to solve the problem.

The method consists in finding $x_{a}(t)$ the solution to the problem that begins in $x(0)=38.6896$ and uses the minimum harvest $h_{min}=10$ until it reaches the value $x^{*}(t)=39.6831$ in the time $t_{min}$.

For this purpose,
$$\left\{
\begin{array}{l}
\dot{x}(t)=0.71 x(t)\left(1-\frac{x(t)}{80.5}\right)-10\\
x(0)=38.6896.\\
\end{array}\right.
$$
Its solution is
\begin{equation}
x_a(t)= \frac{20.9703+62.3013 \, e^{0.388979t}}{1.1523+e^{0.388979t}}.
\label{sola}
\end{equation}
And $t_{min}$ will be such that $x_a (t_{min})= x^{*}(t)=39.6831$. Results $t_{min}=0.232235$. 

Following with the method, it must be found $x_{b}(t)$ solution to the problem that uses the minimum harvest $h_{min}=10$ until it reaches the value $x (T)=x(10)=40.25$.

For this purpose,
$$\left\{
\begin{array}{l}
\dot{x}(t)=0.71 x(t)\left(1-\frac{x(t)}{80.5}\right)-10\\
x(10)=40.25.\\
\end{array}\right.
$$
Its solution is
\begin{equation}
x_b(t)=\frac{889.931+62.3013 \, e^{0.388979t}}{48.9008+e^{0.388979t}}.
\label{solb}
\end{equation}
And $t_{max}$  will be such that $x_b (t_{max})= x^{*}(t)=39.6831$. Results $t_{max}=9.8678$. 

Finally, the optimal solution will be obtained, using (\ref{xestrella}), (\ref{sola}) and (\ref{solb}), times $t_{min}$ and $t_{max}$ and the Nearest Feasible Paths theorem, it is obtained,
\begin{equation}
x^{*}(t)= \left\{ \begin{array}{lcc}
            \displaystyle \frac{20.9703+62.3013 \, e^{0.388979t}}{1.1523+e^{0.388979t}} &   \mathrm{if}  & 0\leq t\leq 0.232235 \\
             \\ 39.6831 &  \mathrm{if} & 0.232235 \leq t\leq 9.8678 \\
             \\ \displaystyle \frac{889.931+62.3013 \, e^{0.388979t}}{48.9008+e^{0.388979t}} &  \mathrm{if}  & 9.8678  \leq t\leq 10=T
             \end{array}
   \right.
\end{equation}
and the optimal harvest obtained with this procedure is
\begin{equation}
h^{*}(t)= \left\{ \begin{array}{lcc}
             10 &   \mathrm{if}  & 0\leq t\leq 0.232235 \\
             \\ 14.2859 &  \mathrm{if} & 0.232235 \leq t\leq 9.8678 \\
             \\ 10 &  \mathrm{if}  & 9.8678  \leq t\leq 10=T.
             \end{array}
   \right.
\end{equation}
Solving the problem in its fractional version as in the previous section, from (\ref{condescalfa}) and (\ref{condeschalfa}),  
is obtained 
\begin{equation}
\begin{array}{ll}
x^{*}_{0.6}(t)&= 40.25\left(1-\frac{0.01}{0.71} (10-t)^{0.4} e^{-0.01\, (10-t)} E_{1,1.4}(0.01\, (10-t))\right).
\end{array}
\end{equation}
While the optimal harvest is
\begin{equation}
h^{*}_{0.6}(t)=0.71 \,x^{*}_{0.6}(t)\left(1-\frac{x^{*}_{0.6}(t)}{80.5}\right)-\,^{C}_{0} D^{0.6}_{t}\left[x^{*}_{0.6}\right](t).
\end{equation}
In the graphic below, the optimal populations corresponding to both cases can be observed.

\begin{figure}[htb]
\centering
\includegraphics[width=85mm]{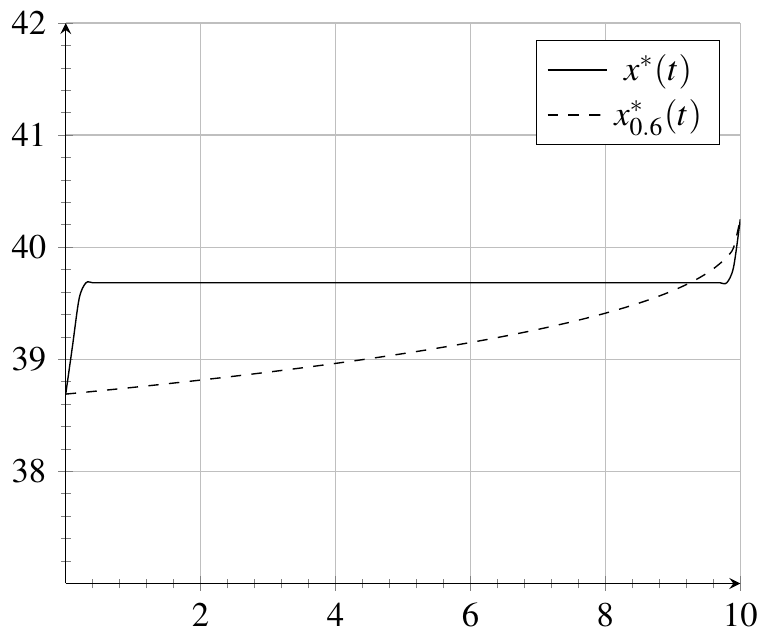}
\caption{Optimal populations of classic and fractional problems.}
\label{Figcosechadescx}
\end{figure}
It can be noticed in Figure \ref{Figcosechadescx} that the optimal solution given by the fractional problem, $x^{*}_{0.6}(t)$, is lower than the optimal population of the classical problem, $x^{*}(t)$, most of the time, it means that the version with $\alpha=0.6$ 
shows a deterioration of the state of the stock with respect to the case $\alpha=1$ which is only recovered at the end by the fact that it has to verify the final condition.

In the following graphics, the optimal harvests are considered. Since $\,^{C}_{0} D^{0.6}_{t}\left[x^{*}_{0.6}\right](t)$ cannot be obtained exactly, we will proceed to use a fractional numerical method of L1 type \cite{BaDiScTr, LiZe}.

A regular partition of $[0,t]$ is considered, as $0=t_0 \leq t_1 \leq$ ... $\leq t_m=t$, of size $\Delta t>0$, to approximate the Caputo derivative as follows
\[\,^{C}_{0} D^{\alpha}_{t}\left[f\right](t_{m}) =\displaystyle \sum_{k=0}^{m-1} b_{m-k-1}(f(t_{k+1}-f(t_k)),\]
where $b_k=\frac{\Delta t^{-\alpha}}{\Gamma(2-\alpha)}\left[(k+1)^{1-\alpha}-k^{1-\alpha}\right]$.

Using this method we have obtained.	

\begin{figure}[htb]	
\centering
\includegraphics[width=85mm]{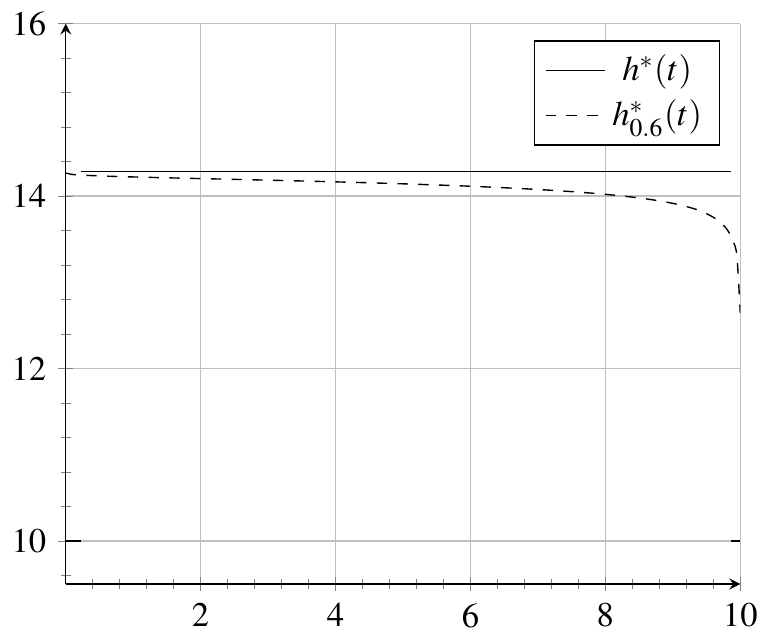}
\caption{Optimal harvests of classical and fractional problems.}
\label{Figcosechadesch}
\end{figure}
It can be noticed in Figure \ref{Figcosechadesch} that although with the fractional problem the extraction of the resource is minor, which is logical because in the fractional model the resource grows more slowly, it only decreases at the end of the interval near $T$ and until the final extraction in $T$ turns out being larger than in the classical case.

This way, it is possible to make a comparison of the profit obtained in each case,
\begin{center}
Classical case profit: $ 134.411$.\\
Fractional case profit:$133.828$.\\
\end{center}
To make an analysis of this, consider the following graphic of the resources evolution without harvest with the given initial condition.

\begin{figure}[htb]
\centering
\includegraphics[width=85mm]{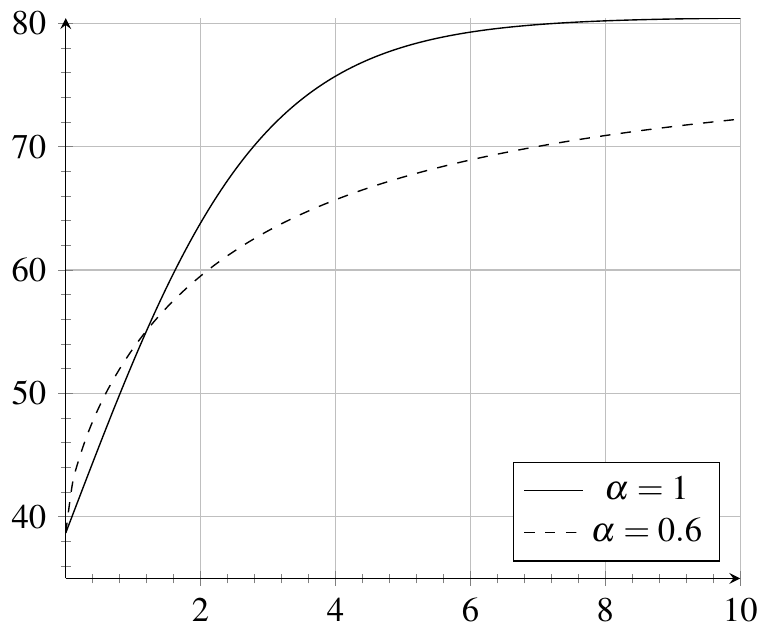}	
\caption{Populations of the classical and fractional problems without harvest.}
\label{Figlog}
\end{figure}
It can be noted that the use of a fractional dynamic equation, which makes the resource grow more slowly as in Figure \ref{Figlog}, does not vary considerably the profit compared to the classical case.

Lastly, it can be stated what happens to the population $x(t)$ if we take the optimal harvest of the classic problem $h_1^{*}$ and consider it in the fractional dynamic equation of the resource.

It must be solved
$$\left\{
\begin{array}{l}
\,^{C}_{0} D^{\alpha}_{t}\left[x\right](t)=0.71\,x(t)\left(1-\frac{x(t)}{80.5}\right)-14.2859\\
x(0)=38.6896.\\
\end{array}\right.
$$
Since this equation has no exact solution, it will be approximated using the Adams fractional method, which consists of using Euler's method to obtain $u_{n+1}^{P}$ (predictor), and the trapezoidal fraction rule to get $u_{n + 1}$ (corrector),
\[ \left\{ \begin{array}{ll}
u_{n+1}^{P}&=\displaystyle \sum_{j=0}^{m-1}\frac{t_{n+1}^{j}}{j!}u_{0}^{j}+\sum_{j=0}^{n}b_{j,n+1}f(t_j,u_j),\\
u_{n+1}&=\displaystyle \sum_{j=0}^{m-1}\frac{t_{n+1}^{j}}{j!}u_{0}^{j}+\sum_{j=0}^{n}a_{j,n+1}f(t_j,u_j)+a_{n+1,n+1}f(t_{n+1},u_{n+1}^{P}).\\
\end{array}\right.\]
For more details refer to \cite{BaDiScTr, LiZe}.

The following result was obtained
\begin{figure}[htb]
\centering
\includegraphics[width=85mm]{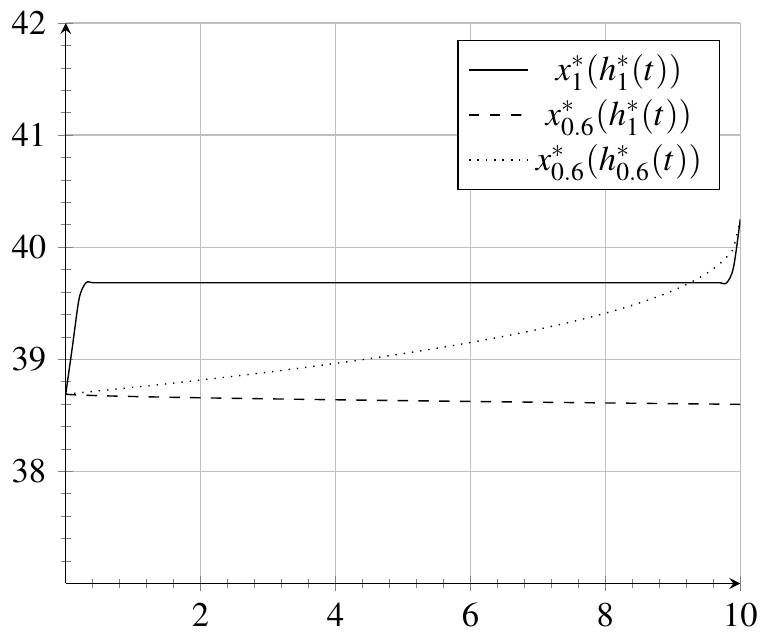}		
\caption{Optimal populations of the classical and fractional with different harvest problems.} 
\label{Figxconh1}
\end{figure}

In the Figure \ref{Figxconh1} it can be seen that the population obtained from the fractional dynamic equation corresponding to taking the optimal harvest of the classical control problem is lower than the obtained by taking the optimal harvest of the fractional problem, as we expected.

Furthermore, if we assume that the "true evolution" of the resource is considering $\alpha = 0.6$ in the dynamic equation and that the harvesting agency considers the dynamic equation of $\alpha = 1$ to be erroneous, the loss is not very great because the agency will use $h^{*}_1(t)$ and its profit will be $133.828$ and not $ 134.411$ as would have been estimated. Also note that the difference between the profits obtained could be more significant if the instantaneous profit function, which in this case is only the harvest, was a little more complex and depended on the stock too as in Remark 1.

\section{Conclusions}
\label{sec:conclusions}
In this article, we have studied a fractional control problem that models the maximization of the profits obtained by exploiting a certain resource. An explanation of the proposed model has been made. Due to the singularity of the problem, different resolution techniques have been developed, both for the classic case and the fractional case. Although we have seen the need of a non-existent fractional Nearest Feasible Paths theorem, we have been able to make a comparison between the classical and fractional results for a certain value of the fractional order. It is also observed that the order of time fractional derivative significantly affects the population growth. Hence, we conclude that fractional derivatives may be more suitable for modeling the evolution of natural resources that naturally have a resilience problem. As a future investigation, it is proposed that the extension of the fractional Nearest Feasible Paths theorem should be explored and the optimal control problem should be extended for more complex instantaneous profits functions. 

\section*{Acknowledgments}
 This work was partially supported by Universidad Nacional de Rosario through the projects ING568 ``Problemas de Control \'Optimo Fraccionario''. The first author was also supported by CONICET through a PhD fellowship.

\end{document}